\documentclass[11pt]{article}
\usepackage{graphicx}
\usepackage{psfrag}
\usepackage{bbm}
\usepackage{amsmath}
\usepackage{amsthm}
\usepackage{url}
\newtheorem*{theo}{Theorem}
\newtheorem{lem}{Lemma}

\title{NP-completeness of Partial Chirotope Extendibility}
\author{Patrick Baier}

\begin{document}
\maketitle

\begin{abstract}
In the monograph ``Axioms and Hulls'' (1992)~\cite{Kn92} Donald Knuth studies some axiomatizations of geometric
situations. The structures described by one of the axiom systems are called CC-systems.
Knuth proves that it is NP-complete to decide, whether a partially defined CC-system
can be extended to a complete CC-system. The aim of this note is to show that Knuth's proof of 
this result also implies that it is NP-complete to decide the extendability
of partially defined chirotopes.
\end{abstract}

\noindent
Together with Stefan Felsner I've been interested in the question, whether partially given
signotopes are extendable (cf.~\cite{FW01} for signotopes). 
Through a talk of Bernd G\"artner in Berlin we learned about a
proof for the NP-completeness of chirotope extendibility by Falk Tschirschnitz~\cite{Tsch03, Tsch01}.
Remarks of Ileana Streinu and Stefan Felsner brought my
attention to an NP-completeness proof of Knuth~\cite{Kn92}.
In fact, the work of Knuth implies that it is NP-complete to decide the extendability
of partially defined chirotopes. The aim of this note is to make the connections clear.

Knuth defines \textbf{CC-systems} as a boolean function of all ordered triples of a given finite groundset (whose elements are usually referred to as points) following five Axioms, where an expression of the form $pqr$ means that the value associated with $(p,q,r)$ is {\itshape true} (p.~3-4)%%
\footnote{Page numbers always refer to Knuth's monograph~\cite{Kn92}.}:\\[5pt]
\textbf{Axiom 1:} $pqr \Longrightarrow qrp$\hfill(cyclic symmetry)\\
\textbf{Axiom 2:} $pqr \Longrightarrow \neg~prq$\hfill(antisymmetry)\\
\textbf{Axiom 3:} $pqr \vee prq$\hfill(nondegeneracy)\\
\textbf{Axiom 4:} $tqr \wedge ptr \wedge pqt \Longrightarrow pqr$\hfill(interiority)\\
\textbf{Axiom 5:} $tsp \wedge tsq \wedge tsr \wedge tpq \wedge tqr \Longrightarrow tpr$\hfill(transitivity)\\
\textbf{Axiom 5':} $tps \wedge tqs \wedge trs \wedge tpq \wedge tqr \Longrightarrow tpr$\hfill(dual transitivity)\\
(Each of these axioms is to be read with an implied quantification ``for all pairwise distinct points'').
Axiom 5' is a dual version of Axiom 5, Knuth shows that they imply each other (p.~5).

A \textbf{pre-CC-system} is analogously defined by Axioms 1, 2, 3 and 5 (p.~11). As to show the equivalence of Axioms 5 and 5' we do only need Axioms 1, 2 and 3 we can equally define pre-CC-systems by using Axiom 5' instead of Axiom 5.
%In pre-CC-systems (and thus in CC-systems) Axiom 5 can be replaced by a dual version:\\
%Axiom 5': $stp \wedge stq \wedge str \wedge tpq \wedge tqr \Longrightarrow tpr$.\\
\textbf{Uniform chirotopes} are defined by Axioms 1, 2, 3 and the so-called \textbf{Grass\-mann-Pl\"ucker-Relations (GPR)} stating that for pairwise distinct points $t, p, q, r, s$ the set 
\begin{eqnarray} 
\label{GPR}
\{(tpq \wedge trs) \vee (tqp \wedge tsr);\notag\\
 (trp \wedge tqs) \vee (tpr \wedge tsq);\\
 (tps \wedge tqr) \vee (tsp \wedge trq)\}\notag
\end{eqnarray}
 contains both values {\itshape true} and {\itshape false}.

Given a pre-CC-system, for an arbitrary point $t$ Knuth considers the {\itshape associated tournament}, which is a complete directed graph on all points except for $t$ defined by the relation $p\rightarrow q$ iff $tpq$ (p.~7). A tournament is called {\itshape vortex-free} iff the graphs shown beneath are not contained as subgraphs (p.~11-12).
\begin{figure}[h]
\centering
\includegraphics[scale=0.4]{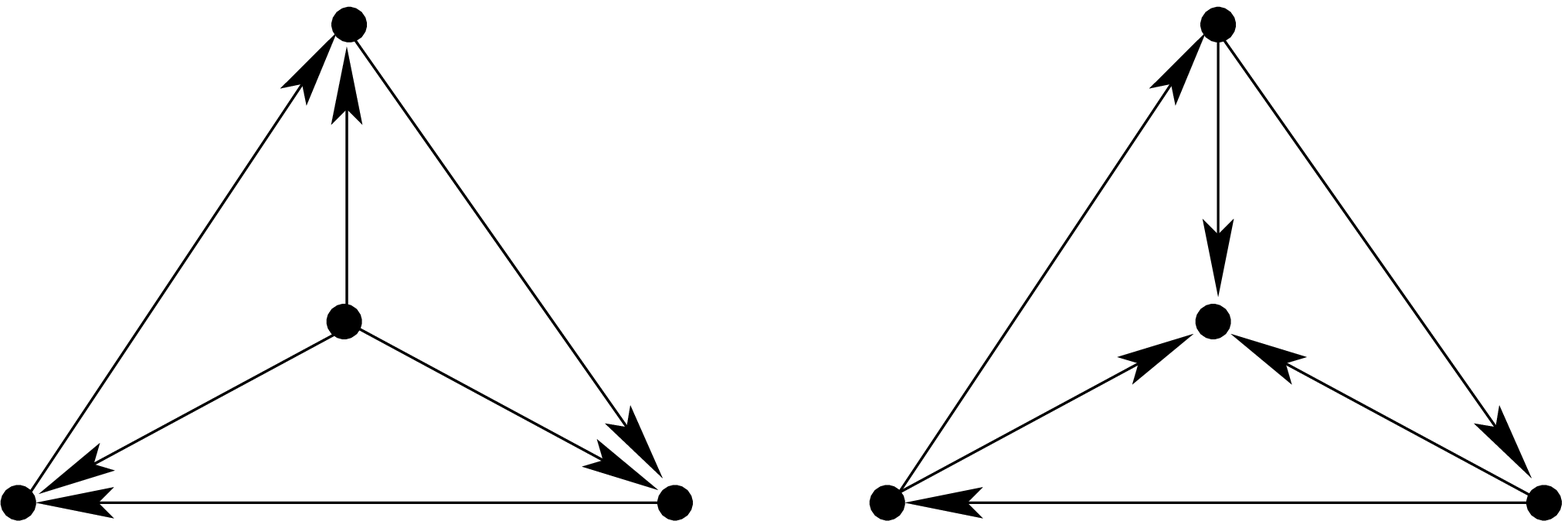}\\
\end{figure}
Knuth now gives the statement of the following
\begin{lem}
\label{vf} 
A pre-CC-system is characterized by the fact that each associated tournament is vortex-free.
\end{lem}
\begin{proof}
This simply follows from the observation that the forbidden subgraphs correspond exactly with the structures excluded by Axioms 5 (left picture) and 5' (right picture), where the tournament associated with $t$ is to be considered, $p, q, r$ correspond to the triangle points in clockwise order and $s$ corresponds to the point in the middle of the respective figure.\footnote{According to the duality of Axioms 5 and 5' it would in fact be enough to exclude one of the pictured subgraphs. The existence of one forbidden structure implies the existence of the other (in a tournament associated with a different point). We explicitly exclude both structures because it simplifies the proof of the following lemma.}
\end{proof}
%Knuth considers another forbidden subgraph corresponding to a dual version of Axiom 5, which follows from Axioms 1, 2, 3, 5 and together with Axioms 1, 2, 3 implies Axiom 5 (p.~5). In his notation what we call {\itshape vortex-free} would be {\itshape out-vortex-free}. For simplicity's sake, we won't consider this dual approach.}\\

At this point we can observe the following simple equivalence:
\begin{lem} The uniform chirotopes are exactly the pre-CC-systems.
\label{equivalence}
\begin{proof} Given Axioms 1, 2, 3 we have to show the equivalence of Axiom 5 (or equivalenty 5 and 5') and the GPR:

%Therefore, it's worth expressing those statements in a less formal way. Let us consider an arbitrary order of the five elements in use, say $pqrst$. 
First let Axiom 5 be violated. That means we have $tsp \wedge tsq \wedge tsr \wedge tpq \wedge tqr \wedge trp$ for some five pairwise distinct points. One can easily check that the GPR aren't fulfilled (all three values in (\ref{GPR}) become {\itshape false}).

For the other direction suppose the GPR are not fulfilled. Using Lemma~\ref{vf} we want to show the existence of a forbidden subgraph. We have two consider two main cases: 
(a) All three values of (\ref{GPR}) are {\itshape true}. This happens if, e.~g.~in either of the three expressions the left part of the $\vee$-expression is {\itshape true}, i.~e.~$tpq \wedge trs \wedge trp \wedge tqs \wedge tps \wedge tqr$. Regarding the tournament associated with $t$, this gives us a subgraph of the form of the right picture, where $p, q, r$ are the triangle points in clockwise order and $s$ is the middle point. Instead of examining all eight cases, where all three values in the GPR become true, we can observe that we can generate all these possibilities from the above mentioned example by sequentially performing an operation, which consists of inverting two disjoint arcs. Looking at the effect of this operation on the tournament associated with $t$, we observe that a forbidden subgraph is always transformed in another forbidden subgraph. Thus starting with the one special case examined above the existence of forbidden subgraphs in all eight cases follows. (b) All three values of (\ref{GPR}) are {\itshape false}. This is the case, if e.~g.~$tpq \wedge tsr \wedge trp \wedge tsq \wedge tsp \wedge tqr$. As tournament associated with $t$ we now get the subgraph shown to the left, where again $p, q, r$ are the triangle points in clockwise order and $s$ is the middle point. In the same way as in part (a) we can avoid examining each of the eight configurations that can occur in this case.
\end{proof}
\end{lem}

The NP-completeness proof for the extendibility of CC-systems consists mainly of two steps (p.~19-23). In the first step it is shown to be NP-complete to decide, whether a directed graph can be extended to a vortex-free tournament by reducing 3SAT (via some steps) to this problem (p.~19-22). As a corollary (p.~22-23) NP-completeness for CC-systems is obtained by using a construction yielding a CC-system (on $n$ points) starting with an arbitrary vortex-free tournament (on $n-1$ points) such that this vortex-free tournament is associated with the additional point. Thus, starting with a set of defined triples (obeying Axioms 1 and 2), all of which have a point $t$ in common, following the same rule as used for the definition of associated tournaments ($p\rightarrow q$ iff $tpq$) we obtain a directed graph on all points except for $t$. This graph can be extended to a vortex-free tournament iff the given set of defined triples can be extended to a CC-system.  
With these results of Knuth it is very easy to show the following:

\begin{theo} The problem of deciding, whether a boolean function defined on a subset of the set of all triples of a given groundset can be extended to a chirotope (or equivalently pre-CC-system), is NP-complete. 
\end{theo}

\begin{proof}
We use the reduction of Knuth. If the given 3SAT instance is not
satisfiable, then there is no completion of the digraph to a
vortex-free tournament and, hence, no completion of the given triples
to a pre-CC-system. If, however, the given 3SAT instance is
satisfiable, then there is a completion to a vortex-free
tournament. Following Knuth we can construct a CC-system, this tournament being associated with one of its nodes. As every CC-system is a pre-CC-system we have found an extension
of the given partial pre-CC-system.
\end{proof}

\end{document}